\magnification=\magstep1
\font\scaps=cmcsc10
\font\msam=msam10
\font\smsam=msam7
\newfam\maths
\textfont\maths=\msam
\scriptfont\maths=\smsam
\def\hexnumber#1{\ifcase#1 0\or 1\or 2\or 3\or 4\or 5\or 6\or 7\or 8\or
9\or A\or B\or C\or D\or E\or F\fi}
\edef\myf{\hexnumber\maths}
\mathchardef\rst="2\myf16
\mathchardef\ltrg="2\myf45
\mathchardef\concat="2\myf61
\def\frac#1/#2{\leavevmode\kern.1em
    \raise.5em\hbox{\the\scriptfont0 #1}\kern-.1em
    /\kern-.15em\lower.25ex\hbox{\the\scriptfont0 #2}}

\def\Lim{{\rm Lim\,}}
\def\Lim{{\rm Lim\,}}
\def\Lim{{\rm Lim\,}}
\def\a{\alpha}
\def\a{\alpha}
\def\bin{{}^{\omega>}2}
\def\bin{{}^{\omega>}2}
\def\bin{{}^{\omega>}2}
\def\can{{}^\omega2}
\def\can{{}^\omega2}
\def\can{{}^\omega2}
\def\cn{\mathbin{\raise2.5pt\hbox{$\concat$}}}
\def\cn{\mathbin{\raise2.5pt\hbox{$\concat$}}}
\def\cn{\mathbin{\raise2.5pt\hbox{$\concat$}}}
\def\d{\delta}
\def\enc#1{\langle#1\rangle}
\def\enc#1{\langle#1\rangle}
\def\enc#1{\langle#1\rangle}
\def\ep{\epsilon}
\def\ep{\epsilon}
\def\eqdf{=^{\rm df}}
\def\eqdf{=^{\rm df}}
\def\eqdf{=^{\rm df}}
\def\fin{{\rm fin}}
\def\fin{{\rm fin}}
\def\fin{{\rm fin}}
\def\g{\gamma}
\def\hide#1{}
\def\len{{\rm length\,}}
\def\len{{\rm length\,}}
\def\len{{\rm length\,}}
\def\nbar{\overline n}
\def\nbar{\overline n}
\def\nbar{\overline n}
\def\nl{\hfill\break}
\def\nl{\hfill\break}
\def\nl{\hfill\break}
\def\om{\omega}
\def\om{\omega}
\def\om{\omega}
\def\seq#1{{}^{#1}2}
\def\seq#1{{}^{#1}2}
\def\seq#1{{}^{#1}2}
\def\si{\sigma}
\def\si{\sigma}
\def\si{\sigma}
\def\su{\subseteq}
\def\su{\subseteq}
\def\su{\subseteq}
\def\tit#1{\noindent{\bf#1}}
\def\tit#1{\noindent{\bf#1}}
\def\tit#1{\noindent{\bf#1}}
\def\z{\zeta}
\def\z{\zeta}
\mathchardef\cge="2\myf3E
\mathchardef\cge="2\myf3E
\mathchardef\cle="2\myf34
\mathchardef\cle="2\myf34

\centerline{{\bf Every null-additive set is meager-additive}\footnote\dag
{ \sevenrm Publ.~No.~445.
First version written in April 91.
Partially supported by the Basic Research Fundations of the
Israel Academy of Sciences and the Edmund Landau Center for research
in Mathematical Analysis, supported by the Minerva Foundation (Germany)
\nl
\null\hskip\parindent
We thank Azriel Levy for rewriting the paper.
The reader should thank him for including more proofs and making the
paper self contained.}}
\smallskip

\centerline{\scaps Saharon Shelah}
\medskip

\tit{\S1. The basic definitions and the main theorem.}
\smallskip

\tit{1. Definition.} (1) We define addition on $\can$ as addition modulo
2 on each component, i.e., if $x,y,z\in\can$ and $x+y=z$ then for every
$n$ we have $z(n)=x(n)+y(n)\pmod 2$.\nl
(2) For $A,B\subseteq\can$ and $x\in\can$ we set $x+A\eqdf\{x+y:y\in A\}$,
and we define $A+B$ similarly.\nl
(3) We denote the Lebesgue measure on $\can$ with $\mu$.
We say that $X\subseteq\can$ is {\it null-additive} if for every
$A\su\can$ which is null, i.e. $\mu(A)=0$, \ $X+A$ is null too.\nl
(4) We say that $X\su\can$ is {\it meager-additive} if for every
$A\su\can$ which is meager also $X+A$ is meager.
\smallskip

\tit{2. Theorem.} Every null-additive set is meager-additive.
\smallskip

\tit{3. Outline and discussion.} Theorem 2 answers a question of Palikowsi.
It will be proved in \S2.
In \S3 we shall present direct characterizations of the null-additive
sets, and in \S4 we shall do the same for the meager-additive sets.

It is obvious that every countable set is both null-additive and
meager-additive.
Are there uncountable null-additive sets, and even null-additive sets of
cardinality $2^{\aleph_0}$?
It will be shown in \S5 that if the continuum hypothesis holds then
there is such a set.
Haim Judah has shown that there is a model of ZFC in which all the
null-additive sets are countable, but there are in it
uncountable meager-additive sets.
This is the model obtained by adding to L more than $\aleph_1$ Cohen
reals.
In this model the Borel conjecture holds, and therefore every
null-additive set is strongly meager and hence countable.
On the other hand, in this model the uncountable set of all constructible
reals is meager-additive.
\medbreak

\tit{\S2 The proof of Theorem 2.}
\smallskip

\tit{4. Notation.} (1) we shall use variables as follow:
$i,j,k,l,m,n$ for natural numbers, $f,g,h$ for functions from $\om$
to $\om$, \ $\eta, \zeta, \nu, \sigma, \tau$ for finite sequences of 0's and
1's, $x,y,z$ for members \hbox{of $\can$}, \ $A,B,X,Y$ for subsets of $\can$,
and $S,T$ for trees.\nl
(2) $\bin=\bigcup_{n<\omega}{}^n2$.
We shall denote subsets of $\bin$ with $U,V$.
For $\eta\in\bin$, $U\su\bin$ and $x\in\can$ we shall write $\eta+x$ for
$\eta+x\rst\len(\eta)$, and $U+x$ for $\{\eta+x:\eta\in U\}$.\nl
(3) For $\eta,\nu\in\bin$ we write $\eta\ltrg\nu$ if $\nu$ is an
extension of $\eta$.\nl
(4) A tree is a nonempty  subset of $\bin$ such that
\item{} (a) If $\eta\ltrg\nu$ and $\nu\in T$ then also $\eta\in T$, and
\item{} (b) If $\eta\in T$ and $n>\len(\eta)$ then there is a $\nu$ of
length $n$ such that $\eta\ltrg\nu$ and $\nu\in T$.

\noindent
(5) For a tree $T$ \ $\Lim T=\{x\in\can:\hbox{\rm for every }n<\om
\enspace x\rst n\in T\}$.\nl
(6) A tree $T$ is said to be {\it nowhere dense} if for every  $\eta\in T$
there is a $\tau\in\bin$ such that $\eta\ltrg\tau$ and $\tau\notin T$.
A set $B\su\can$ is said to be nowhere dense if $B\su\Lim T$ for some
nowhere dense tree $T$.\nl
(7) For every $x,y\in\can$ we write $x\equiv y$ if $x(n)=y(n)$ for all
but finitely many $n<\omega$.
For $A\su\can$ \
$A^\fin\eqdf\{y\in\can:y\equiv x \;\hbox{\rm for some } x\in A\}$.\nl
(8) $U^{[\nu]}\eqdf\{\tau\in U:\tau\ltrg\nu\>{\rm or}\>\nu\ltrg\tau\}$
\enspace(read: $U$ through $\nu$).\nl
(9) $U^{\enc\nu}\eqdf\{\tau\in\bin:\nu\cn\tau\in U\}$
\enspace (read: $U$ above $\nu$),
and for $\eta\in\bin$ \nl
$\eta^{\enc k}\eqdf\langle\eta(k+i):i<\len\eta-k\rangle$.\nl
(10) For $\nu,\eta\in\bin\cup\can$ we write $\nu\sim_n\eta$ if
$\len(\nu)=\len(\eta)$ and $\nu(i)=\eta(i)$ for every $n\le i<\len(\nu)$.
For $S\su\bin\cup\can$ we define
$S^{\sim n}=\{\nu:\nu\sim_n\eta\hbox{ for some }\eta\in S\}$.
\smallbreak

\tit{5. Outline of the proof.}
Let $X\su\can$ be null-additive.
It clearly suffices to prove that for every $A\su \can$ which is nowhere
dense $X+A$ is meager.
Given a nowhere dense tree $S$ we shall prove in Lemma 6 a condition
which is sufficient for a tree $T$ to be such that $T+S$ is nowhere
dense.
Then we shall split $X$ to a union $X=\bigcup_{i=1}^\infty X_i$ such that
for each $i$ \ $X_i\su\Lim T_i$ where $T_i$ is a tree which satisfies
that condition.
Thus for a nowhere dense $S$ each set $X_i+\Lim S\su\Lim(T_i+S)$ is
nowhere dense, hence $X+\Lim S\su\bigcup_{i=1}^\infty\Lim(T_i+S)$ is
meager.
\smallbreak

\tit{6. Lemma.} Let $T$ be a tree such that\nl
(a) $T$ is nowhere dense.\nl
(b) $f=f_T$ is the function from $\om$ to $\om$ given by\nl
$f(n)=\min\{m:\hbox{for every }\eta\in\seq n\hbox{ there is a
}\tau\in\seq m\hbox{ such that }\eta\ltrg\tau\hbox{ and }\tau\notin
T\}$.\nl
Thus for every sequence $\eta$ of length $n$ there is a witness of
length $\le f(n)$ that $T$ is nowhere dense.
Obviously for every $n<\om$ \ $f(n)>n$, and if $n<m$ then
$f(n)\le f(m)$.

Let $g$ be a function from $\om$ to $\om$ and $\nbar=\enc{n_i:i<\om}$,
$\nbar'=\enc{n'_i:i<\om}$ increasing sequences of natural numbers such
that\nl
(c) $f^{g(i)}(n_i)\le n'_i<n_{i+1}$ for every $i<\om$, where $f^m$
denotes the $m$-th iteration of $f$.\nl
Then for every tree $S$ which satisfies\nl
(d) $S$ is of width $(\nbar',g)$, i.e., for every $i<\om$ \
$|\seq{n'_i}\cap S|\le g(i)$,\nl
$T+S$ is nowhere dense.
\smallbreak

\noindent
{\it Proof.} Let $\eta\in\seq{n_i}$.
We shall show the existence of an $\eta'\in\seq{n'_i}$ such that
$\eta\ltrg\eta'$ and $\eta'\notin T+S$.

By (c) there is a sequence $m_0,\dots,m_{g(i)}$ such that $m_0=n_i$, \
$f(m_k)\le m_{k+1}$ for $0\le k\le g(i)$ and $m_{g(i)}=n'_i$.
Let $\enc{\tau_k:k<k_i}$ enumerate the set $\seq{n'_i}\cap S$.
$k_i\le g(i)$ by~(d).
We define $\eta_k\in\seq{m_k}$ for $0\le k\le k_i$ by recursion as
follows.
$\eta_0=\eta$.
Given $\eta_k\in\seq{m_k}$, for $k<k_i$, we shall define
$\eta_{k+1}\in\seq{m_{k+1}}$ so that for no extension
$\eta'\in\seq{n'_i}$ of $\eta_{k+1}$ we shall have $\eta'+\tau_k\in T$.
We have $\eta_k+\tau_k\rst m_k\in\seq{m_k}$ and by the definition of $f$
and by the choice of the
$m_k$'s $\eta_k+\tau_k\rst m_k$ has an extension $\nu\in\seq{m_{k+1}}$
such that $\nu\notin T$.
If we take $\eta_{k+1}=\nu+\tau_k\rst m_{k+1}$ then
$\eta_k+\tau_k\rst m_k\ltrg\nu$ implies $\eta_k\ltrg\eta_{k+1}$ ,
\ $\eta_{k+1}\in\seq{m_{k+1}}$
and $\eta_{k+1}+\tau_k\rst m_{k+1}=\nu\notin T$, and therefore
for every $\eta'\in\seq{n'_i}$ such that $\eta_k\ltrg\eta'$ we have
$\eta'+\tau_k\notin T$.
Let $\eta'=\eta_{k_i}$, and assume that $\eta'\in T+S$.
Then, for some $k<k_i\le g(i)$ \ $\eta'+\tau_k\in T$, contradicting our
choice of $\eta_{k+1}=\eta'\rst m_{k+1}$.
Thus $\eta'\notin T+S$.
\smallbreak

\tit{7. Lemma.} If $S,T_i,\>i\in\om$ are trees and
$\Lim S\su\bigcup_{i\in\om}\Lim T_i$ then for some $\eta\in S$ and
$j\in\om$ \ $S^{[\eta]}\su T_j$.
\smallskip

\noindent
{\it Proof.} Suppose that this is not the case, i.e., for every
$\eta\in S$ and $i<\om$ there is a $\zeta$ such that
$\zeta\in S^{[\eta]}$ and $\zeta\notin T_i$.
Once there is such a $\zeta$ we can assume that
$\eta\ltrg\zeta$ and $\len\zeta>\len\eta$.
We define now, by induction on $i$, \ $\eta_i$ and $k_i$ so that
$k_i=\len\eta_i$, \ $k_0=0$, \ $\eta_0=\enc{}$, \ $\eta_i\ltrg\eta_{i+1}$,
\ $k_i<k_{i+1}$, \ $\eta_{i+1}\in S$ and $\eta_{i+1}\notin T_i$.
Let $y=\bigcup_{i\in\om}\eta_i$.
Since $\eta_i\in S$ for every $i\in\om$ \ $y\in\Lim S\su\bigcup_{i\om}\Lim T_i$,
hence for some $j\in\om$ \ $y\in\Lim T_j$.
However, $y\rst k_{j+1}=\eta_{j+1}\notin T_j$, contradicting $y\in\Lim T_j$.
\smallbreak

\tit{8. Lemma.} Let $S$ and $T$ be trees such that
$\Lim S\su(\Lim T)^\fin$.
Then there are $k<\om$, \ $\eta,\nu\in\seq k$, \ $\eta\in S$ such that
$S^{\enc\eta}\su T^{\enc\nu}$.
\smallskip

\noindent
{\it Proof.}
For $n<\om$, \ $\si_1,\si_2\in\seq n$ and $\si_2\in T$ we define\nl
$T_{\si_1,\si_2}\eqdf\{\tau:\tau\ltrg\si_1\}\cup
\{\si_1\cn\tau:\si_2\cn\tau\in T\}$  (This is the tree $T^{[\si_2]}$
with ``$\,\si_2$ replaced by $\si_1$'').
Clearly
$$(\Lim T)^\fin=
\bigcup_{n<\om,\, \si_1,\si_2\in\seq n,\, \si_2\in T}\Lim T_{\si_1,\si_2}
\leqno(1)$$
Since there are only countably many  $T_{\si_1,\si_2}$'s in (1)
there are by Lemma 7 a $\zeta\in S$ and $j<\om$ such that
$S^{[\zeta]}\su T_{\si_1,\si_2}$.
Clearly there is an $\eta$ with $\zeta\ltrg\eta$ and a $\nu$ with
$\len\nu=\len\eta$ such that $S^{\enc\eta}\su T^{\enc\nu}$.
(If $\zeta\ltrg\si_1$ then $\eta=\si_1$ and $\nu=\si_2$, else
$\si_1\ltrg\zeta$ and then $\eta=\zeta$ and
$\nu=\si_2\cn\zeta\rst[\len\zeta,\len\si_2)\>$).
\smallbreak

\tit{9. Lemma.} Let $X$ be a null-additive set.
Let $T$ be a tree such that $\mu(\Lim T)>0$\hide{, and let $M$ be an unbounded
subset of $\om$}.
There is a tree $T^*$ such that $\mu(\Lim T^*)>0$, for every $\eta\in T^*$
\ $\mu(\Lim(T^{*[\eta]}))>0$, and
$((\can\setminus(\Lim T)^\fin)+X)\cap\Lim T^*=\emptyset$, and then
$X=\bigcup_{\eta\in T^*,\>\len\zeta=\len\eta\hide{\in M}}Y_{\eta,\zeta}$ where
$Y_{\eta,\zeta}=\{x\in X:\zeta\cn x^{\enc{\len\zeta}}+T^{*[\eta]}\su T\}$.
\smallskip

\noindent
{\it Proof.} Since $\mu(\Lim T)>0$ then, as easily seen,
$\mu((\Lim T)^\fin)=1$, hence\nl
$\mu(\can\setminus(\Lim T)^\fin)=0$.
Since $X$ is null-additive also $\mu(X+(\can\setminus(\Lim T)^\fin))=0$.
Hence there is a tree $T^*$ such that $\mu(\Lim T^*)>0$ and
$(X+(\can\setminus(\Lim T)^\fin))\cap\Lim T^*=\emptyset$.
Without loss of generality we can assume that $T^*$ has been pruned so
that for $\eta\in T^*$ \ $\mu(\Lim T^{*[\eta]})>0$.

Let $x\in X$ then
$\can\setminus(x+(\Lim T)^\fin)=x+(\can\setminus(\Lim T)^\fin)\su
X+(\can\setminus(\Lim T)^\fin)$.
Hence $(\can\setminus(x+(\Lim T)^\fin))\cap\Lim T^*\su
(X+(\can\setminus(\Lim T)^\fin)\cap\Lim T^*=\emptyset$, i.e.,\break
$\Lim T^*\su x+(\Lim T)^\fin$, and therefore
$\Lim(x+T^*)=x+\Lim T^*\su(\Lim T)^\fin$.
By Lemma 8 there are $\eta\in T^*$ and $\nu\in\seq{\len\eta}$ such that
$x^{\enc{\len\eta}}+T^{*\enc{\eta}}\su T^{\enc\nu}$.
\hide{If $\len\eta\notin M$ we can extend $\eta$ and $\nu$ to sequences
whose length in $M$ and which still satisfy the same inclusion.}
Let $\zeta=\eta+\nu$, then $\zeta+\eta=\nu$ and therefore
$\zeta\cn x^{\enc{\len\eta}}+T^{*[\eta]}\su T^{[\nu]}\su T$, hence
$x\in Y_{\eta,\zeta}$.
\smallbreak

\tit{10. Lemma.} Let $X$ be null-additive, and let
$\overline n=\enc{n_i:i<\om}$, \ $\overline n'=\enc{n'_i:i<\om}$
be such that for every $i<\om$ \ $n_i<n'_i$ and
$n'_i+i\cdot 2^{n_i'}\le n_{i+1}$,
then we can represent $X$ as $\bigcup_{m<\om}X_m$
such that for each $m$, for some real $a_m \in (0, 1)$ and $S_m$
of width $(\overline n',g_{a_m})$ we have $X_m\su\Lim(S_m)$,
where for every real $a\in(0,1)$ \
$g_a$ is the function on $\om$ given by $g_a(0)=1$, \
$g_a(i)={\rm max}(1,{\rm int}(\log_2(a)/\log_2(1-2^{-i})))$,
where for a real $d$ \
${\rm int}(d)$ is the integral part of $d$.
\smallskip

\noindent
{\it Proof.} Since $n'_i+i\cdot2^{n'_i}<n_{i+1}$
we can fix for each $0<i<\om$ a sequence
$\enc{u_{i,\tau}:\tau\in\seq{n'_i}}$
of pairwise disjoint subsets of the interval $[n'_i,n_{i+1})$
having $i$ members each.
Let $B\su\can$ be given by
$$ B=\{y\in\can:(\forall j>0)(\exists k\in u_{j,y\rst n'_j})\,y(k)=1\}.$$
$B$ is clearly a closed subset of $\can$ hence for
$T=\{y\rst n:y\in B\wedge n\in\om\}$ \ $B=\Lim(T)$.

The properties of $T$ in which we are interested are\nl
(B0) $T\supseteq\seq{n_1}$.\nl
(B1) For each $\eta\in T\cap\seq{n'_i}$ \
$\left|T^{[\eta]}\cap\seq{n_{i+1}}\right|=2^{(n_{i+1}-n'_i)}(1-2^{-i})$.\nl
(B2) If $\eta,\nu_0,\dots,\nu_{k-1}\in\seq{n'_i}$, \
$\nu^+_0,\dots,\nu^+_{k-1}\in\seq{n_{i+1}}$, \
$\eta+\nu_l\in T$, \ $\nu_l\ltrg\nu^+_l$ for $l<k$ 
and $\nu_0,\dots , \nu_{k-1}$ is with no repetitions
then
$\left|\{\eta^+:\eta\ltrg\eta^+\in\seq{n_{i+1}},\>
(\forall l<k)(\eta^++\nu_l^+\in T)\}\right|\le
2^{n_{i+1}-n'_i}\left(1-2^{-i}\right)^k$.\nl
(B3) For every  $\eta \in \seq{n'_i}$ we have : $\eta \rst n_i \in T$
implies $\eta \in T$.

These properties can be established by an obvious counting argument.

By (B0), (B1) and (B3) we have\nl
$$\eqalign{\mu(\Lim T)&=\mu\left(\bigcap_{i=1}^\infty\{x\in\can:x\rst
n_i\in T\}\right)\cr
&=\mu\left(\{x\in\can:x\rst n_1\in T\}\right)\cdot
\prod_{i=1}^\infty{\mu(\{x\in\can:x\rst n_{i+1}\in T\})\over
\mu(\{x\in\can:x\rst n_i\in T\})}\cr
&=1\cdot\prod_{i=1}^\infty{|T\cap\seq{n_{i+1}}|/2^{n_{i+1}}\over
|T\cap\seq{n_{i+1}}|/2^{n_i}}=\prod_{i=1}^\infty(1-2^{-i})>0\cr}$$

For the $T$ which we constructed let $T^*$ and $Y_{\eta,\zeta}$ be as in
Lemma 9.
For $\rho\in\seq{\len\eta}$ let
$Y_{\eta,\zeta,\rho}=\{y\in Y_{\eta,\zeta}:y\rst\len\eta=\rho\}$.
Clearly
$$X=\bigcup_{\eta\in T^*,\>\len\eta=\len\zeta=\len\rho}Y_{\eta,\zeta,\rho}
\leqno(2)$$
Since there are only countably many $Y_{\eta,\zeta,\rho}$'s they can be
taken to be the $X_m$'s we are looking for, provided we show that every such
$Y_{\eta,\zeta,\rho}$ is a subset of $\Lim(S)$ for some tree $S$ of width
$\enc{\overline n',g_a}$ for some real $0<a<1$.
We shall see that this is indeed the case if we take
$S=\{y\rst m: y\in Y_{\eta,\zeta,\rho},\,m<\om\}$ and $a=\mu(T^{*[\eta]})$.
$a>0$ by what we assumed about $T^*$.
As, obviously, $Y_{\eta,\zeta,\rho}\su\Lim(S)$ all we have to do is to show
that $S$ is of width $\enc{\overline n',g_a}$.

We can choose a set $W \subseteq S \cap \seq{n_{j+1}}$  such that
the function mapping $\eta \in W$ to $\eta \rst n'_j$ is one to one
and onto $S \cap \seq{n'_j}$

We fix now $\eta,\zeta,\rho$ and denote $Y_{\eta,\zeta,\rho}$ by $Y$ and
the length of $\eta,\zeta,\rho$ by $n$.
Let $z\in\can$ be such that $z\rst n=\zeta+\rho$ and $z(i)=0$ for $i\ge n$.
Then for every $y$ such that $y\rst n=\rho$ we have $y+z=\zeta\cn y^{\enc n}$.
Therefore, by the definition of $Y$ we have
$$Y=\{y\in\can:y\rst n=\rho,\;(\zeta\cn y^{\enc n})+T^{*[\eta]}\su T\}
=\{y\in\can:y\rst n=\rho,\;y+z+T^{*[\eta]}\su T\}\leqno(3)$$
for every $y\in Y$
there is a unique $\tau\in W$ such that $\tau\rst n'_j=y\rst n'_j$ \
($\tau$ may be $y\rst n_{j+1}$).
Clearly $|W|=|S\cap\seq{n'_j}|$ and we denote $|W|$ with $s$,
so it suffices to prove $s\le g_a(j)$.
If $n'_j\le n$ then the only member of $S\cap\seq{n'_j}$ is $\rho\rst n'_j$
hence $s=1$, so $s\le g_a(j)$.
We shall now deal with the case where $n'_j>n$.
Let $\tau_0,\dots\tau_{s-1}$ be the members of $W$.
For $m<s$ \ $\tau_m=y\rst n_{j+1}$ for some $y\in Y$, hence, by
(3), $\tau_m+z+T^{*[\eta]}\su T$ and therefore
$(z+T^{*[\eta]})\cap\seq{n_{j+1}}\su\tau_m+T$.
Since this holds for every $\tau\in W$ we have
$$ z+T^{*[\eta]}\cap\seq{n_{j+1}}\su\bigcap_{m<s}\tau_m+T\leqno(4)$$

Let us find out the size of $\bigcap_{m<s}(\tau_m +T)$.
Let $\sigma\in\seq{n'_j}$, and we shall ask how many members $\tau$ of
$\bigcap_{m<s}(\tau_m + T)$ extend $\sigma$.
Now $\tau\in\tau_m+T$ for each $m<s$
iff $\tau+\tau_m\in T$ for each $m<s$.
If for some $m<s$ \ $\sigma+\tau_m\rst n'_j\notin T$ then
also $\tau+\tau_m\notin T$, hence $\sigma$ has no extension in
$\bigcap_{m<s}(\tau_m+T)$.
If for every $m<s$ \ $\sigma+\tau_m\rst n'_j\in T$ then
by (B2) (where $\eta=\si$, \ $\nu_m=\tau_m\rst n'_j$ and
$\nu_m^+=\tau_m$), since $\tau_m\rst n'_j\neq\tau_l\rst n'_j$ for $m\neq l$,
the number of $\tau$'s such that $\si\ltrg\tau\in\seq{n_{j+1}}$
and $\tau+\tau_m\in T$ for every $m<s$ is $2^{n_{j+1}-n'_j}(1-2^{-j})^s$.
Since there are $2^{n'_j}$ different $\si$'s in $\seq{n'_j}$ we have
$$\left|\bigcap_{m<s}(\tau_m+T)\right|\le
2^{n_{j+1}}\cdot(1-2^{-j})^s.\leqno(5)$$

On the other hand, since $\mu(T^{*[\eta]})= a$ \
$T^{*[\eta]}\cap\seq{n_{j+1}}$ has at least $a\cdot 2^{n_{j+1}}$
members, and so has $z+T^{*[\eta]}\cap\seq{n_{j+1}}$.
Comparing (4) with (5) we get
$a\cdot2^{n_{j+1}}\le2^{n_{j+1}}(1-2^{-j})^s$, i.e.,
$a\le(1-2^{-j})^s$, \ $\log_2(a)\le s\cdot\log_2(1-2^{-j})$, \
$s\le\log_2(a)/\log_2(1-2^{-j})$.
\medbreak

\tit{11. Proof of Theorem 2.} Let $X$ be null-additive.
As mentioned in {\bf 5} it suffices to
show that for every nowhere dense tree $T$ \ $X+\Lim(T)$
is meager.
Let $f = f_T$ as in Lemma 6.
Define by recursion $n_0=0$, \ $n'_i = f^{g_{1/(i + 1)}(i)}(n_i)+1$ and
$n_{i+1}= n'_i+i\cdot2^{n'_i}+1$.
By Lemma 10 $X\su\bigcup_{m<\om}\Lim(S_m)$,
where for some $a_m \in (0,1)$ \ $S_m$ is of width
$\enc{\overline{n}',g_{a_m}}$, hence it
suffices to show that if $S$ is of width $\enc{\overline{n}',g_a}$
for some $a\in(0,1)$ then $\Lim(S)+\Lim(T)=\Lim(S+T)$ is meager.
Let $j$ be such that ${1\over j+1}\le a$ and let
$\eta_1,\dots,\eta_k$ be all the members of $S$ of length $n'_j$.
Then $S=\bigcup_{l=1}^k S^{[\eta_l]}$ and
$\Lim S=\bigcup_{l=1}^k\Lim(S^{[\eta_l]})$.
Therefore it suffices to prove that for $1\le l\le k$ \
$\Lim (S_l)+\Lim(T)$ is meager and this follows once we show that
$S_l+T$ is nowhere dense.
To prove this we show that the requirements of Lemma 6 hold here for
$S_l, T$.
(a) and (b) hold by our choice of $T$ and $f$.
Let $g$ be defined by $g(i)=1$ for $i<j$ and $g(i)=g_a(i)$ for $i\ge j$.
Now we shall see that (c) holds.
For $i<j$ we have
$n'_i=f^{g_{1/(i+1)}(i)}(n_i)+1\ge f(n_i)+1=f^{g(i)}(n_i)+1$,
since $f(n)\ge n$ for every $n$, and for $i\ge j$ we have
$n'_i=f^{g_{1/(i+1)}(i)}(n_i)+1\ge f^{g_a(i)}(n_i)+1=f^{g(i)}(n_i)+1$,
since $a\ge{1\over j+1}\ge{1\over i+1}$ and $g_a(i)$ is a
decreasing function of $a$.
Thus for every $i<\om$ \
$f^{g(i)}(n_i)\le f^{g_{1/(i+1)}(i)}(n_i)\le n'_i$.
(d) of Lemma 6 holds since for $i<j$ \ $|\seq{n'_i}\cap S_l|=1=g(i)$
and for $i\ge j$ \
$|\seq{n'_i}\cap S_l|\le |\seq{n'_i}\cap S|\le g_a(i)=g(i)$.

%
%
\tit{\S3 Characterization of the null-additive sets}
\smallskip

\tit{12. Definition.} By a {\it corset\/} we mean a non decreasing function
$f$ from $\om$ to
$\om\setminus\{0\}$ which converges to infinity (i.e., for every $n<\om$
\ $f(m)>n$ for all sufficiently large $m$).
For a corset $f$, we say that a tree $T$ is {\it of width $f$} if for
every $n<\om$ \ $|T\cap\seq n|\le f(n)$; and we say that $T$ is
{\it almost of width} $f$ if $|T\cap\seq n|\le f(n)$ for all
sufficiently large $n$.
\smallbreak

\tit{13. Theorem.}
For every $X\su\can$ the following conditions are equivalent:\nl
a. $X$ is null-additive.\nl
b. For every corset $f$ there is a tree $S$ of width $f$ such that
$X\su\Lim(S)^{\fin}$.\nl
c. For every corset $f$ there are trees $S_m$, \ $m<\om$, which are
almost of width $f$ such that $X\su\bigcup_{m<\om}\Lim(S_m)^{\fin}$.\nl
d. For every corset $f$ there are  trees $S_m$, \ $m<\om$,
of width $f$ such that $X\su\bigcup_{m<\om}\Lim(S_m)$.
\smallbreak

\tit{Proof.}
(b)$\rightarrow$(c) is obvious.
\smallbreak

\noindent
(c)$\rightarrow$(d).
Let $S$ be a tree almost of width $f$.
Then for some $k$ we have $|T\cap\seq n|\le f(n)$ for all $n\ge k$.
By~(1) of Lemma~8 $\Lim(S)^\fin=
\bigcup_{\sigma_1,\sigma_2\in\seq k,\>\sigma_2\in S}
\Lim(S_{\sigma_1,\sigma_2})$.
Each $S_{\sigma_1,\sigma_2}$ is of width $f$ since for $n\le k$ we have
$|S_{\sigma_1,\sigma_2}\cap\seq n|=1$ and for $n>k$ we have
$|S_{\sigma_1,\sigma_2}\cap\seq n|\le|S\cap\seq n|\le f(n)$.
Therefore, if $X\su\bigcup_{m<\om}\Lim(S_n)^\fin$ as in (c) then each
$S_m$ can be replaced by countably many $S_{\sigma_1,\sigma_2}$'s and
(d) holds.
\smallbreak

\noindent
(d)$\rightarrow$(b).
Let $f$ be a corset.
We can easily define by recursion a sequence $0=n_0<n_1<\dots$ of
natural numbers and a corset $f^*$ such that for all $j<\om$ and
$m\ge n_{j+1}$ we have $(j+1)\cdot2^{n_j}\cdot f^*(m)\le f(m)$.

For a given corset $f$, if $X$ satisfies (d) let $S^*_m$, \ $m<\om$,
be as in (d) for the corset $f^*$.
We construct now a set $S\su\bin$ by defining $S\cap\seq m$ by recursion
on $m$.
$S\cap\seq0=\{\enc{}\}$.
For $n_i<m\le n_{i+1}$ let
$$S\cap\seq m=\{\eta\in\seq m:\eta\rst n_i\in S\cap\seq{n_i}\hbox{ and }
\eta\in S_j^{*\sim n_j} \hbox{ for some } j < i \vee j=0 \}.$$
$S$ can be easily seen to be a tree, and clearly
$\Lim(S)^\fin\supseteq\bigcup_{n<\om}\Lim(S^*_n)\supseteq X$.
For $m\le n_1$ easily , for  $n_i\le m<n_{i+1} ,i\ge 1$ we have\nl
$|S\cap\seq m|\le\sum_{j\le i}\left|S_j^{*\sim n_j}\cap\seq m\right|
=\sum_{j\le i}2^{n_j}|S_j^*\cap\seq m|\le
(i+1)\cdot 2^{n_j} \cdot f^*(m) \le f(m)$,\nl
thus $S$ is of width $f$.
\smallbreak

\noindent
(d)$\rightarrow$(a).
Assume now that (d) holds for $X$,
and let $A\su\can$, \ $\mu(A)=0$; we shall prove that $\mu(X+A)=0$.
First we shall mention two lemmas of measure theory the proof of which
is left to the reader.

\tit{Lemma A.} For every tree $T$ with
$\mu(\Lim(T))=a>0$ and $\ep>0$ there is an $N\in\om$ such that for every
$n\ge N$ there is a $t\subseteq\seq n\cap T$ such that
$|t|\ge2^n(a-\ep)$ and for each $\eta\in t$ \
$\mu(\Lim(T^{[\eta]})>2^{-n}(1-\ep)$.

Using Lemma A one can prove

\tit{Lemma B.} For every
tree $T$ with $\mu(\Lim(T))>0$, every $\ep>0$ and every sequence
$\enc{\ep_i:0<i<\om}$ of positive reals there is a subtree $T'$ of $T$ and
an increasing sequence $\enc{n_i:i<\om}$ of natural numbers such that
$n_0=0$, \ $\mu(\Lim(T'))>\mu(\Lim(T))-\ep$ and
$$\hbox{for $i>0$ and every }\eta\in\seq{n_i}\cap T'\quad
\mu(\Lim(T'{}^{[\eta]}))>2^{-n_i}(1-\ep_i)\leqno(6)$$

By basic mesure theory $\mu(A^{\fin})=0$ so
there is a tree $T$ such that $\mu(\Lim(T))>0$ and
$\Lim(T) \cap A^\fin =\emptyset$ hence
$\Lim(T)^\fin\cap A=\emptyset$.
Given $\ep<\mu(\Lim(T))$ and $\enc{\ep_i:i<\om}$ as in Lemma B we obtain a
subtree $T'$ of $T$ as in that lemma with $\mu(\Lim(T'))>0$.
The union of sufficiently many ``finite translates'' of $T'$, i.e.,
trees $T'_{\si_1,\si_2}$ as in (1) of Lemma 8 is a tree $T''$ satisfying
(6) with $\mu(\Lim(T''))\ge{1\over2}$. \
$\Lim(T'')^\fin=\Lim(T')^\fin\su\Lim(T)^\fin$ and hence
$\Lim(T'')\cap A\su\Lim(T)^\fin\cap A=\emptyset$.
We take now $\ep_i={1\over4(i+1)^3}$ and take $T$ to be $T''$ and
we get $\mu(\Lim(T))\ge{1\over2}$ and
$$\hbox{for $i>0$ and every }\eta\in\seq{n_i}\cap T\quad
\mu(\Lim(T {}^{[\eta]}))>2^{-n_i}(1-{1\over4(i+1)^3})\leqno(7)$$

Let $f$ be the corset given by $f(n)=i+1$ for $n_i\le n<n_{i+1}$.
By (d) there are trees $S_m$ of width $f$ such that
$X\su\bigcup_{m<\om}\Lim(S_m)$.
To show that $\mu(X+A)=0$ it clearly suffices to show that
for every tree $S$ of width $f$ \ $\mu(\Lim(S)+A)=0$.

We define
$$T^*=\{\eta\in\bin:\nu+\eta\in T
\hbox{ for every $\nu\in S$ of the same length as $\eta$}\}$$
We do not show that $T^*$ is a tree but obviously if
$\zeta\ltrg\eta\in T^*$ then $\zeta\in T^*$, thus $\Lim(T^*)$ is
defined. 
If $\mu(\Lim(T^*))>0$ then, by a well-known property of the measure,
$\mu(\Lim(T^*)^\fin)=1$, hence in order to prove $\mu(\Lim(S)+A)=0$
it suffices to prove $(\Lim(S)+A)\cap\Lim(T^*)^\fin=\emptyset$.
Assume $y\in(\Lim(S)+A)\cap\Lim(T^*)^\fin$.
Since $y\in\Lim(T^*)^\fin$ there is a $y'\in\can$ such that $y'(n)=y(n)$
for all sufficiently big $n$'s and $y'\in\Lim(T^*)$.
Since $y\in\Lim(S)+A$ there is an $x\in\Lim(S)$ such that $y+x\in A$,
hence $y+x\notin\Lim(T)^\fin$, hence $y'+x\notin\Lim(T)$.
Therefore, for some $n$ \ $y'\rst n+x\rst n\notin T$, hence, by the
definition of $T^*$, \ $y'\rst n\notin T^*$ contradicting
$y'\in\Lim(T^*)$.

We still have to prove that $\mu(\Lim(T^*))>0$.
We shall prove, by induction on $i$, that
$$n_i\le n\le n_{i+1}\rightarrow
|(T\setminus T^*)\cap\seq
n|\le2^n\cdot\sum_{j<i}{1\over4(j+1)^2}.\leqno(8)$$
Once we establish (8) we notice that since\nl
$\Lim(T)\setminus\Lim(T^*)=
\bigcup_{n<\om}\Lim(T)\setminus\{x\in\can:x\rst n\in T^*\}$, and
the set\nl
$\Lim(T)\setminus\{x\in\can:x\rst n\in T^*\}$ is increasing with
$n$ hence $\mu(\Lim(T)\setminus\Lim(T^*))$\nl
$=\lim_{n\rightarrow\infty}\mu(\Lim(T)\setminus\{x\in\can:x\rst n\in T^*\})\le
\lim_{n\rightarrow\infty}2^{-n}|(T\setminus T^*)\cap\seq n|$\nl
$\le\lim_{n\rightarrow\infty}\sum_{j=0}^n{1\over4(j+1)^2}=
\sum_{j=0}^\infty{1\over4(j+1)^2}={\pi^2\over24}<{1\over2}$
and since $\mu(\Lim(T))\ge{1\over2}$ \nl
$\mu(\Lim(T^*))>0$.

To prove (8), assume now $n_i\le n\le n_{i+1}$.
By the definition of $T^*$\nl
$(T\setminus T^*)\cap\seq n=
\{\eta\in T\cap\seq n:(\exists\rho\in S\cap\seq n)\rho+\eta\notin T\}$\nl
$=\{\eta\in T\cap\seq n:
(\exists\rho\in S\cap\seq n)(\eta\rst n_i+\rho\rst n_i\notin T)\}\cup$\nl
\null\hfill$\bigcup_{\rho\in S\cap\seq n}\{\eta\in T\cap\seq n:
\eta\rst n_i+\rho\rst n_i\in T\,\wedge\,\eta+\rho\notin T\}$\break
$\su\{\eta\in\seq n:\eta\rst n_i\in T\setminus T^*\}\cup
\bigcup_{\rho\in S\cap\seq n}\{\eta\in\seq n:
\eta+\rho\in\{\si\in\seq n:\si\rst n_i\in T\,\wedge\,\si\notin T\}\}$.\nl
Therefore $|(T\setminus T^*)\cap\seq n|\le
2^{n-n_i}|(T\setminus T^*)\cap\seq{n_i}|+|S\cap\seq n|
|\{\si\in\seq n:\si\rst n_i\in T\wedge\si\notin T\}|$.
For $i>0$ we have, by the induction hypothesis
$|T\setminus T^*\cap\seq{n_i}|\le2^{n_i}\sum_{j<i}{1\over4(j+1)^2}$.
For $i=0$ we have $(T\setminus T^*)\cap\seq{n_i}=\emptyset$ since
$n_0=0$ and $\emptyset\in T^*$.
$|S\cap\seq n|\le f(n)=i$ and
$|\{\si\in\seq n:\si\rst n_i\in T\wedge\si\notin T\}|\le
{2^n\over4(i+1)^3}$, by (7).
Thus\nl
$|(T\setminus T^*)\cap\seq n|\le
2^{n-n_i}\cdot2^{n_i}\sum_{j<i}{1\over4(j+1)^2}+(i+1)\cdot{2^n\over4(i+1)^3}
\le2^n\sum_{j<i+1}{1\over4(j+1)^2}$ which is what we had to show.
\smallbreak

\noindent
(a)$\rightarrow$(c).
Most of the proof follows that of Lemma 10.
We need also the following Lemma 14, which will be proved later.
Let $f$ be a corset.
\smallbreak

\tit{14. Lemma.}
There is an infinite sequence $0=n_0<n_1<n_2<\dots$ and a tree $T$ such that
for every $i\in\om$ \ $f(n_{i+1})>(i+1)\cdot2^{i+1}+1$ and\nl
(B1) For each $\eta\in T\cap\seq{n_i}$ we have
$|T^{[\eta]}\cap\seq{n_{i+1}}|=2^{(n_{i+1}-n_i)}\cdot(1-2^{-(i+1)})$.\nl
(B2) If $\eta,\nu_0,\dots,\nu_{k-1}\in\seq{n_i}$, \
$\nu^+_0,\dots,\nu^+_{k-1}\in\seq{n_{i+1}}$, \
$\nu^+_j\neq\nu^+_l$ for $j<l<k$, \
$\eta+\nu_l\in T$, \ $\nu_l\ltrg\nu^+_l$ for $l<k$ then\nl
$\left|\{\eta^+:\eta\ltrg\eta^+\in\seq{n_{i+1}},\>
(\forall l<k)(\eta^++\nu_l^+\in T)\}\right|\le
2^{n_{i+1}-n_i}\left(1-2^{-(i+1)}\right)^{k-1}$.
\smallbreak

Let $\enc{n_i:i\in\om}$ and $T$ be as in Lemma 14.
As in the proof of Lemma 10 we get $\mu(\Lim T)>0$.
Let $T^*$ and $Y_{\eta,\zeta}$ be as in Lemma 9 and let
$Y_{\eta,\zeta,\rho}\,$, \ $S$ and $z$ be as in the proof of Lemma~10.
All~we have to do is to show that $S$ is almost of width $f$.
Let us fix $\eta$, \ $\zeta$ and $\rho$.
We shall now see that
$$\eqalign{&\hbox{If $\eta'\in T^{*[\eta]}\cap\seq{n_i}$ then }\cr
&|\{\eta^+:\eta'\ltrg\eta^+\in T^*\cap\seq{n_{i+1}}\}|/
2^{(n_{i+1}-n_i)}
\le(1-2^{-(i+1)})^{|S\cap\seq{n_i}|-1}}\leqno(9)$$
Let $\eta^+\in T^{*[\eta]}\cap\seq{n_{i+1}}$, then, by the
definition of $S$ (see (3)), if $\rho^+\in S\cap\seq{n_{i+1}}$ then
$\rho^++\eta^++z\in T$.
Thus\nl
$\{\eta^+:\eta'\ltrg\eta^+\in T^*\cap\seq{n_{i+1}}\}\su
\{\eta^+:\eta'\ltrg\eta^+\in\seq{n_{i+1}},\>
(\forall\rho^+\in S)\rho^++\eta^++z\in T\}$.\nl
Let us take in (B2) $\eta=\eta'$, \ $k=|S\cap\seq{n_i}|$, \
$\{\tau_l:l<k\}=S\cap\seq{n_i}$,
$\{\tau_l^+:l<k\}\su S\cap\seq{n_{i+1}}$, and for $l<k$ \
$\tau^+_l\rst n_i=\tau_l$, \ $\nu_l=\tau_l+z$, \
$\nu_l^+=\tau_l^++z$, hence
$\nu_l=\nu_l^+\rst n_i$ for $l<k$.
Since for $l<k$ \ $\nu_l^++z=\tau_l^+\in S\cap\seq{n_{i+1}}$ we have\nl
$\{\eta^+:\eta'\ltrg\eta^+\in\seq{n_{i+1}},\>
(\forall\rho^+\in S)(\rho^++\eta^++z\in T\}$\nl
\null\hfill $\su\{\eta^+:\eta'\ltrg\eta^+\in\seq{n_{i+1}},\>
(\forall l<k)(\nu_l^++\eta^+\in T)\}$,\break
therefore by (B2)\nl
$|\{\eta^+:\eta'\ltrg\eta^+\in\seq{n_{i+1}},\>
(\forall\rho^+\in S)(\rho^++\eta^++z\in T\}|\le
2^{n_{i+1}-n_i}(1-2^{-(i+1)})^{|S\cap\seq{n_i}|-1}$,
which establishes (9).

(9) tells us how $T^*$ grows from the level $n_i$ to the level $n_{i+1}$
and therefore\nl
$|T^*\cap\seq{n_i}|\cdot2^{-n_i}\le
\prod_{j<i}(1-2^{-(j+1)})^{|S\cap\seq{n_j}|-1}$.\nl
Let $c_0=\mu(\Lim T^*)$.
We know that $c_0>0$ and we can assume $c_0<1$.
Then\nl
$-\infty<\log c_0\le\log(|T^*\cap\seq{n_i}|\cdot2^{-n_i})
\le\sum_{j<i}\bigl(\log(1-2^{-(j+1)})\cdot(|S\cap\seq{n_j}|-1)\bigr)$.
Since $\log(1-x)\le-{1\over2}x$ we get
$\sum_{j<i}2^{-(j+2)}\cdot(|S\cap\seq{n_{i+1}}|-1)\le\log{1\over c_0}$.
We shall denote $4\log{1\over c_0}$ by $c$, so
$\sum_{j<i}2^{-j}\cdot(|S\cap\seq{n_j}|-1)\le c$, and for every $j$ \
$2^{-j}(|S\cap\seq{n_j}|-1)\le c$, hence
$|S\cap\seq{n_j}|\le c\cdot2^j+1$.
For $j>c$ we have, by our choice of the $n_i$'s,\nl
$f(n_j)>j\cdot2^j+1>c\cdot2^j+1\ge|S\cap\seq{n_j}|$, hence $S$ is almost
of width $f$.
\smallbreak

Lemma 14 follows immediately from the following Lemma.
\smallskip

\tit{15. Lemma.}
For every $n\in\om$ and $0<p<1$ there is an $N>n$ such that for every
$n'\ge N$ and $t\su\seq n$ there is a $t'\su\seq{n'}$ which satisfies
the following (i)--(iii).\nl
(i) For each $\z\in t'$ \ $\z\rst n\in t$.\nl
(ii) For each $\eta\in t$ \ $|t'{}^{[\eta]}|\ge2^{n'-n}\cdot p$.\nl
(iii) If $0<k\le2^n$ \ $\eta,\nu_0,\dots,\nu_{k-1}\in\seq n$, \
$\nu_0^+,\dots,\nu_{k-1}^+\in\seq{n'}$, \ $\nu^+_j\neq\nu^+_l$ for
$j<l<k$, \ $\eta+\nu_l\in t$, \ $\nu_l=\nu^+_l\rst n$ for $l<k$ then
$|\{\eta^+:\eta\ltrg\eta^+\in\seq{n'},
(\forall l<k)\>\eta^++\nu_l^+\in t'\}|\le2^{n'-n}p^{k-1}$.
\smallskip

\tit{Proof.} We shall prove the lemma by the probabilistic method.
Let $n'>n$ and let $A=\{\eta^+\in\seq{n'}:\eta^+\rst n\in t\}$.
We construct a subset $A^*$ of $A$ as follows.
We take a coin which yields heads with probability $p$.
For each $\eta^+\in A$ we toss this coin and we put $\eta^+$ in $A^*$
iff the coin shows heads.
We shall see that if we take $t'=A^*$ then, for sufficiently large $n'$,
the probability that (ii) holds has a positive lower bound which does not
depend on $n'$ while the probability that (iii) holds is arbitrarily close
to $1$.
Hence there is an $N$ and a $t'$ as claimed by the lemma.
We prove first two lemmas.
\smallbreak

\tit{Lemma 16.}
For $k,\eta,\nu_0,\dots,\nu_{k-1},\nu_0^+,\dots,\nu^+_{k-1}$ as im
Lemma~15 there are reals $c_1,c_2>0$ which depend only on $p$, $n$ and
$k$ such that
$${\rm Pr}\,\Bigl(|\{\eta^+:\eta\ltrg\eta^+\in\seq{n'},\;
\bigwedge_{l<k}\eta^++\nu^+_l\in A*\}|\ge p^{k-1}2^{n'-n}\Bigr)
<c_1e^{-c_2\cdot2^{n'}}$$

\tit{Proof.} We denote $2^{n'-n}$ with $m$.
We set $(\seq{n'})^{[\eta]}=\{\eta_j^+:j<m\}$.
Let $G$ be the graph on $m$ given by
$$iGj\enspace{\rm iff}\enspace\{\eta^+_i+\nu^+_l:l<k\}\cap
\{\eta^+_j+\nu^+_l:l<k\}\neq\emptyset$$
Obviously each $i<m$
has at most $k^2$ neighbors in $G$ hence, by
a well known theorem, $m$ can be decomposed into $k^2+1$ pairwise
disjoint sets $B_0,\dots,B_{k^2}$ such that for every $i\le k^2$ if
$j,l\in B_i$ and $j\neq l$ then $jGl$ does not hold.
Let $d<{1\over2}{\rm
min}\,\{p^{l-1}-p^l:l\le2^n\}={1\over2}p^{2^n-1}(1-p)>0$.
$$\eqalign{&
{\rm Pr}\,\Bigl(|j<m:\bigwedge_{l<k}\eta^+_j+\nu_l^+\in A^*\}|
\ge m\cdot p^{k-1}\Bigr)\cr
&\le{\rm Pr}\,\Bigl(|j<m:\bigwedge_{l<k}\eta^+_j+\nu_l^+\in A^*\}|>
m(p^k+d)\Bigr)\qquad\qquad\hbox{since $p^k+d<p^{k-1}$}
\cr}\leqno(10)$$
Assume that
$$\eqalign{&\hbox{for every }i\le k^2\hbox{ such that }
|B_i|\ge{dm\over2k^2+2}\cr
&\hbox{ we have }|\{j\in B_i:\bigwedge_{l<k}\eta^+_j+\nu_l^+\in A^*\}
\le|B_i|(p^k+{d\over2})\cr}\leqno(11)$$
then
$$\{j<m:\bigwedge_{l<k}\eta^+_j+\nu^+_l\in A^*\}\su
\bigcup_{i\le k^2,\>|B_i|\ge{dm\over2k^2+2}}\{j\in B_i:
\bigwedge_{l<k}\eta_j^++\nu_l^+\in A^*\}
\bigcup_{i\le k^2,\>|B_i|<{dm\over2k^2+2}}B_i$$
hence
$$\eqalign{&|j<m:\bigwedge_{l<k}\eta_j^++\nu_l^+\in A^*\}|\cr
&\le\sum_{i\le k^2,\>|B_i|\ge{dm\over2k^2+2}}
|j\in B_i:\bigwedge_{l<k}\eta_j^++\nu_l^+\in A^*\}|
+\sum_{i\le k^2,\>|B_i|<{dm\over2k^2+2}}|B_i|\cr
&\le\sum_{i\le k^2,\>|B_i|\ge{dm\over2k^2+2}}
|B_i|(p^k+{d\over2})
+\sum_{i\le k^2,\>|B_i|<{dm\over2k^2+2}}|B_i|,\qquad\qquad\qquad
\hbox{by (11)}\cr
&\le m(p^k+{d\over2})+(k^2+1){dm\over2k^2+2}=m(p^k+d)\cr}$$
Therefore the event
$|\{j<m:\bigwedge_{l<k}\eta_j^++\nu_l^+\in A^*\}|>m(p^k+d)$ is
incompatible with (11), so we continue the inequality (10) by
$$\eqalign{&\le{\rm Pr}\>\Bigl(\bigvee_{i\le k^2,\>|B_i|\ge{dm\over2k^2+2}}
\bigl(|\{j\in B_i:\bigwedge_{l<k}\eta_j^++\nu_l^+\in A^*\}|
>|B_i|(p^k+{d\over2})\bigr)\Bigr)\cr
&\le\sum_{i\le k^2,\>|B_i|\ge{dm\over2k^2+2}}{\rm Pr}\>
\Bigl(|\{j\in B_i:\bigwedge_{l<k}\eta_j^++\nu_l^+\in A^*\}|
>|B_i|(p^k+{d\over2})\Bigr)\cr}\leqno(12)$$
For a fixed $j<m$ the events $\eta_j^++\nu_l^+\in A^*$ for different $l$'s
are independent hence
${\rm Pr}\>\Bigl(\bigwedge_{l<k}\eta_j^++\nu_l^+\in A^*\Bigr)=p^k$.
For a fixed $i$ the events
$\bigwedge_{l<k}\eta_j^++\nu_l^+\in A^*$ for different $j$'s in $B_i$
are independent since, by the definition of the $B_i$'s, if
$j_1,j_2\in B_i$ and $j_1\neq j_2$ then
$\eta_{j_1}^++\nu_{l_1}^+\neq\eta_{j_2}^++\nu_{l_2}^+$.
We have here $|B_i|$ independent events, each with probability $p^k$.
By a formula of probability theory (see, e.~g., the formula
${\rm Pr}\,[X>a]<e^{-2a^2/n}$ in Spencer [2], p.~29)
$${\rm Pr}\>\Bigl(\{j\in B_i:\bigwedge_{k<l}\eta_j^++\nu_l^+\in A^*\}|
>|B_i|p^k+\ep\Bigr)<e^{-{2\ep^2\over|B_i|}}$$
and taking $\ep={1\over2}|B_i|d$ we get
$${\rm Pr}\>\Bigl(\{j\in B_i:\bigwedge_{k<l}\eta_j^++\nu_l^+\in A^*\}|
>|B_i|(p^k+{d\over2})\Bigr)<e^{-{d^2|B_i|\over2}}$$
Continuing (12) we get
$$\le\sum_{i\le k^2,\>|B_i|\ge{dm\over2k^2+2}}e^{-{d^2|B_i|\over2}}
\le\sum_{i\le k^2,\>|B_i|\ge{dm\over2k^2+2}}e^{-{d^2\over2}{dm\over2k^2+2}}
\le(k^2+1)e^{-{d^32^{n'-n}\over4k^2+4}}$$
Combining this with the inequalities (10) and (12) we get
$$\eqalign{&{\rm Pr}\,\Bigl(|\{\eta^+:\eta\ltrg\eta^+\in\seq{n'},\;
\bigwedge_{l<k}\eta^++\nu^+_l\in A^*\}|\ge p^{k-1}2^{n'-n}\Bigr)\cr
&<(k^2+1)e^{-{d^32^{n'-n}\over4k^2+4}}
=(k^2+1)e^{-{d^32^{-n}2^{n'}\over4k^2+4}}\cr}$$
Since $d={1\over2}p^{2^n-1}(1-p)$ this proves Lemma 16.
\smallbreak

\tit{17. Lemma.}
There are $c_3,c_4$ which depend only on $p$ and $n$ such that
$$\eqalign{&{\rm Pr}\Bigl(
\bigvee_{k,\eta,\nu_0,\dots,\nu_{k-1},\nu_0^+,\dots,\nu_{k-1}^+}\cr
&\hbox to2cm{}|\{\eta^+:\eta\ltrg\eta^+\in\seq{n'},\;
(\forall l<k)\>\eta^++\nu_l^+\in A^*\}|\ge2^{n'-n}p^{k-1}\Bigr)\cr
&\le c_3(2^{n'})^{2^n}e^{-c_42^{n'}}\cr}\leqno(13)$$
where $k,\eta,\nu_0,\dots,\nu_{k-1},\nu_0^+,\dots,\nu_{k-1}^+$ are as in
(iii) of Lemma 15.
\smallskip

\tit{Proof.}
By our requirements on
$k,\eta,\nu_0,\dots,\nu_{k-1},\nu_0^+,\dots,\nu_{k-1}^+$ there are at
most $2^n$ possible $k$'s and $\eta$'s and $(2^{n'})^{2^n}$ sequences
$\enc{\nu_0^+,\dots,\nu_{k-1}^+}$, while $\nu_0,\dots,\nu_{k-1}$ are
determined by $\nu_0^+,\dots,\nu_{k-1}^+$ and $n$.
Therefore we get, by Lemma 16,
$$\eqalign{&{\rm Pr}\Bigl(
\bigvee_{k,\eta,\nu_0,\dots,\nu_{k-1},\nu_0^+,\dots,\nu_{k-1}^+}
\bigl(|\{\eta^+:\eta\ltrg\eta^+\in\seq{n'},\;
(\forall l<k)\>\eta^++\nu_l^+\in A^*\}|\ge2^{n'-n}p^{k-1}\bigr)\Bigr)\cr
&\le\sum_{k,\eta,\nu_0,\dots,\nu_{k-1},\nu_0^+,\dots,\nu_{k-1}^+}
{\rm Pr}\>\Bigl(|\{\eta^+:\eta\ltrg\eta^+\in\seq{n'},\;
(\forall l<k)\>\eta^++\nu_l^+\in A^*\}|\ge2^{n'-n}p^{k-1}\Bigr)\cr
&\le2^n\cdot2^n\cdot(2^{n'})^{(2^n)}\cdot c_1e^{-c_22^{n'}}\cr}$$
\smallbreak

\tit{Proof of Lemma 15 (continued).}
For each $\eta^+\in\seq{n'}$ such that $\eta\ltrg\eta^+$ \
$\eta^+\in A^*$ if the coin shows heads and different tosses are
independent $|A^{*[\eta]}|$ is a binomial random variable with
expectation $2^{n'-n}p$.
By the central limit theorem of probability theory (see, e.~g., Feller
[1, Ch.~7]) the limit, as $n'\rightarrow\infty$, of
${\rm Pr}\>\bigl(|A^{*[\eta]}|\ge2^{n'-n}p\bigr)$ is
$\int_0^\infty{1\over\sqrt{2\pi}}e^{-{x^2\over2}}dx={1\over2}$,
hence there is an $N$ such that for every $n'\ge N$ \
${\rm Pr}\>\bigl(|A^{*[\eta]}|\ge2^{n'-n}p\bigr)\ge{1\over3}$.
For different $\eta\in t$ the random variables $|A^{*[\eta]}|$ are
idependent, hence
$${\rm Pr}\>\bigl(\bigwedge_{\eta\in t}(|A^{*[\eta]}|\ge2^{n'-n}p)\bigr)
\ge{1\over3^{|t|}}\ge{1\over3^{2^n}}\leqno(14)$$
The right-hand side of (13) clearly vanishes as $n\rightarrow\infty$, let
us take N to be such that for $n'\ge N$ the right-hand side of (13) is
$<3^{-2^n}$.
Therefore we have, by (13) and (14),
$$\eqalign{&{\rm Pr}\Bigl(
\bigvee_{k,\eta,\nu_0,\dots,\nu_{k-1},\nu_0^+,\dots,\nu_{k-1}^+}\cr
&\hbox to2cm{}(|\{\eta^+:\eta\ltrg\eta^+\in\seq{n'},\;
(\forall l<k)\>\eta^++\nu_l^+\in A^*\}|<2^{n'-n}p^{k-1})\cr
&\hbox to1cm{}\wedge
\bigwedge_{\eta\in t}(|A^{*[\eta]}|\ge2^{n'-n}p)\Bigr)>0\cr}\leqno(15)$$
By (15) there is a $t'$ as required by the lemma.


\tit{\S4 Characterization of the meager-additive sets}
\smallskip

\tit{18. Theorem.} For every $X\su\can$ the following conditions are
equivalent:\nl
a. $X$ is meager additive.\nl
b. For every sequence $n_0<n_1<n_2<\dots$ of natural numbers there is a
sequence\break
$i_0<i_1<\dots$ of natural numbers and a $y\in\can$ such
that for every $x\in X$ and for every sufficiently big $k<\om$ there is
an $l\in[i_k,i_{k+1})$ such that
$x\rst[n_l,n_{l+1})=y\rst[n_l,n_{l+1})$.
\smallbreak

\tit{Proof.}
Throughout this proof, if $x\in\can\cup\bin$, \ $k,l\in\om$ and $k<l$ 
then $x\rst[k,l)$ will denote the sequence $\xi\in\seq{l-k}$ such that
$\xi(i)=x(k+i)$ for all $i<l-k$.

\noindent
(b)$\rightarrow$(a).
In order to prove (a) it clearly suffices to show that $X+\Lim T$ is
meager for every nowhere dense tree $T$.

For a nowhere dense tree $T$ let $\enc{n_i:i<\om}$ be an ascending
sequence of natural numbers such that $n_0=0$ and for every $i\in\om$
there is a sequence $\nu_i\in\seq{n_{i+1}-n_i}$ such that for every
$\tau\in\seq{n_i}$ \ $\tau\concat\nu_i\notin T$.
Let $\enc{i_j:j<\om}$ and $y$ be as in (b), then, by (b),
$X=\bigcup_{k\in\om}X_k$ where
$X_k=\{x\in X:(\forall m\ge k)(\exists l\in[i_m,i_{m+1}))\>
x\rst[n_l,n_{l+1})=y\rst(n_l,n_{l+1})$.
It clearly suffices to prove that $X_k+\Lim T$ is nowhere dense.

Let $\tau\in\seq{n_{i_m}}$ for some $m\ge k$; we shall show that $\tau$
has an extension which is not in $X_k+\Lim T$.
Let $\nu=\nu_{i_m}\concat\nu_{i_m+1}\concat\dots\concat\nu_{i_{m+1}-1}$
and let $\rho=y\rst[n_{i_m},n_{i_{m+1}})+\nu$.
We show that no extension $z$ of $\tau\concat\rho$ is in $X_k+\Lim T$.
Suppose $\tau\concat\rho\ltrg z\in X_k+\Lim T$ then $z=x+w$, \
$x\in X_k$ \ $w\in\Lim T$.
Therefore $\tau=\tau_1+\tau_2$ and $\rho=\rho_1+\rho_2$ such that
$\tau_1\concat\rho_1\ltrg x$ and $\tau_2\concat\rho_2\ltrg w$, hence
$\tau_2\concat\rho_2\in T$.
Let $\xi\in\seq{n_{i_m}}$ be such that $\xi(j)=0$ for every $j<n_{i_m}$,
and let $\rho'=\xi\concat\rho$, \ $\rho'_1=\xi\concat\rho_1$, \
$\rho'_2=\xi\concat\rho_2$.
Clearly $\rho'=\rho'_1+\rho'_2$.
Since $x\in X_k$ there is, by (b), an $l\in[i_m,i_{m+1})$ such that
$x\rst[n_l,n_{l+1})=y\rst[n_l,n_{l+1})$.
Since $\tau_1\concat\rho_1\ltrg x$ we have
$\rho'_1\rst[n_{i_m},n_{i_m+1})=x\rst[n_{i_m},n_{i_m+1})$
and hence
$\rho_1\rst[n_l,n_{l+1})=x\rst[n_l,n_{l+1})=y\rst[n_l,n_{l+1})$ 
Therefore, by the definition of $\rho$ and~$\nu$ \
$$\eqalign{y\rst[n_l,n_{l+1})+\rho'_2\rst[n_l,n_{l+1})&
=\rho'_1\rst[n_l,n_{l+1})+\rho'_2\rst[n_l,n_{l+1})=
\rho'\rst[n_l,n_{l+1})\cr
&=y\rst[n_l,n_{l+1})+\nu_l\cr}$$
hence $\rho'_2\rst[n_l,n_{l+1})=\nu_l$.
By the definition of $\nu_l$ \ $\tau_2\concat\rho_2\notin T$,
contradicting $\tau_2\concat\rho_2\in T$.
\smallbreak

\noindent
(a)$\rightarrow$(b).
Let $X$ be meager-additive.
Let $\enc{n_i:i<\om}$ be an ascending sequence of natural numbers.
Let $B=\{x\in\can:\forall j(\exists k\in[n_j,n_{j+1}))\,x(k)\neq0\}$
and $T=\{x\rst n:x\in B,\>n\in\om\}$.
Clearly $B=\Lim T$ is nowhere dense, so $X+\Lim T$ is meager, hence
there are nowhere dense trees $S_n$, \ $n\in\om$ such that for every $n$
\ $S_n\su S_{n+1}$ and $X+\Lim T\su\bigcup_{n\in\om}S_n$.
We define now the sequences $\enc{i_l:l<\om}$, which is an ascending
sequence of natural numbers, and $\enc{\nu_l:l<\om}$ by recursion as
follows.
$i_0=0$.
Given $i_l$ let $\nu_l$ and $i_{l+1}$ be such that
$\nu_l\in\seq{n_{i_{l+1}}-n_{i_l}}$ and for every $\rho\in\seq{n_{i_l}}$
\ $\rho\concat\nu_l\notin S_l$; there are such $\nu_l$ and $i_{l+1}$
since $S_l$ is nowhere dense.
Let $y\in\can$ be given by $y\rst[n_{i_l},n_{i_{l+1}})=\nu_l$ for every
$l<\om$.
We shall now prove that $\enc{i_l:l<\om}$ and $y$ are as required by
(b).

Let $x\in X$, so $\Lim(x+T)=x+\Lim T\su X+\Lim
T\su\bigcup_{n\in\om}S_n$.
Therefore, by Lemma~7 (where we take $x+T$ for $S$) there is an
$\eta\in T$ and $n\in\om$ such that \hbox{$x+T^{[\eta]}\su S_n$}.
Let $k$ be such that $k\ge n$ and $i_k\ge\len\eta$.
By $x+T^{[\eta]}\su S_n$ we have
$\hbox{$x\rst n_{i_{k+1}}+(T^{[\eta]}\cap\seq{n_{i_{k+1}}})$}\su S_n\su S_k$.
Thus for every $\rho\in T^{[\eta]}\cap\seq{n_{i_{k+1}}}$ \ \
$x\rst n_{i_{k+1}}+\rho\in S_k$, hence, by the definition of $\nu_k$ and
$y$, \
\def\rnik{\rst[n_{i_k},n_{i_{k+1}})}
$x\rnik+\rho\rnik\neq\nu_k=y\rnik$ and therefore
$x\rnik-y\rnik\neq\rho\rnik$, i.e., \
$x\rnik-y\rnik\notin\{\rho\rnik:\rho\in T^{[\eta]}\}$.
Since $i_k>\len\eta$ this can happen, by the definition of $T$, only if
for some $i_k\le j<i_{k+1}$ \ $x\rst[n_j,n_{j+1})-y\rst[n_j,n_{j+1})$ is
identically zero, and this is what we had to prove.
\bigbreak

\tit{\S5. An uncountable null-additive set.}
\smallskip

\tit{19. Theorem.} If the continuum hypothesis holds then there is an
uncountable null-additive set.
\smallskip

\tit{Proof.} Let $\enc{f_\a:\a<\om_1}$ be a sequence containing all
corsets and let $\enc{T_\a:\a<\omega_1}$ be a sequence containing all
perfect trees.
Let $E$ be the set of all limit ordinals $\d<\om_1$ such
that for every $\a,\beta<\delta$ and $n<\om$ there is a $\g<\d$
such that
$$T_\g\su T_\a,\enskip
  T_\g\cap\seq n=T_\a\cap\seq n\enskip
  \hbox{and for all } m\quad
  |T_\g\cap\seq m|\le{\rm max}\,(|T_\a\cap\seq m|,f_\beta(m))$$
Clearly $E$ is closed.
For every $\a,\beta<\om_1$ there is a perfect tree $T$ such that
$T\su T_\a$, \ $T\cap\seq n=T_\a\cap\seq n$ and for all $m<\om$ \
$|T\cap\seq m|\le{\rm max}\,(|T_\a\cap\seq n|,f_\beta(m))$.
This tree $T$ is $T_\g$ for some $\g<\om_1$.
By a simple closure argument this implies that $E$ is unbounded.

We need now the following lemma which will be proved later.
\smallbreak

\tit{20. Lemma.} There is an increasing and continuous sequence
$\enc{\d_\z:\z<\om_1}$ of ordinals in~$E$ such that for every $\z<\om_1$, \
$k<\om$ and $\a<\d_\z$ there is an ordinal $\g$ which is good for 
$(\z,\a,k)$, where by $\g$ {\it is good for} $(\z,\a,k)$ we mean that
$$\eqalign{&({\rm i})\;\g<\d_{\z+1}\cr
&({\rm ii})\;T_\g\su T_\a,\enskip T_\g\cap\seq k=T_\a\cap\seq k\cr
&({\rm iii})\;\hbox{for all $\xi\le\z$ such that $\d_\xi>\a$ and for
every $\ep<\d_\z$, there is a $\beta<\d_\xi$}\cr
&\hbox{\phantom{$({\rm iii})$}
such that $T_\g\su T_\beta\su T_\a$ and $T_\beta$ is almost of width
$f_\ep$}\cr}\leqno(16)$$

For $\xi<\om_1$ let $\g_\xi$ be
the $\g$ which is good for $(\xi,0,0)$.
We choose \hbox{$\eta_\xi\in\Lim T_{\g_\xi}\setminus\{\eta_\beta:\beta<\xi\}$},
and let $X=\{\eta_\xi:\xi<\om_1\}$.
$X$ is clearly uncountable.
We shall prove that $X$ is null-additive by proving that $X$ satisfies
condition (c) of Theorem~13.
For a given corset $f$ \ $f=f_\ep$ for some $\ep<\om_1$.
Let $\xi<\om_1$ be such that $\d_\xi>\ep$.
Let \hbox{$Z=\{\beta<\d_{\xi+1}:T_\beta \hbox{ is almost of width }f_\ep\}$.}
We shall see that\nl
$X\su\{\eta_\z:\z\le\xi\}\cup\bigcup_{\beta\in Z}\Lim T_\beta$.
Since $Z$ and $\xi$ are countable condition (c) of Theorem~13 holds.

Let $\z>\xi$, it suffices to prove that $\eta_\z\in\Lim T_\beta$ for
some $\beta\in Z$.
$\ep<\d_\xi$ and since $\g_\zeta$ is good for $\a=k=0$ hence there is a
$\beta<\d_\xi$ such that $T_{\g_\z}\su T_\beta$, and $T_\beta$ is of
width $f_\ep$.
Thus $\beta\in Z$ and $\eta_\zeta\in\Lim T_{\g_\z}\su\Lim T_\beta$.
\smallbreak

\tit{Proof of Lemma 20.} We define $\enc{\d_\z:\z<\om_1}$ as follows.
$\d_0$ is the least member of $E$.
For a limit ordinal $\z$  \ $\d_\z=\bigcup_{\xi<\z}\d_\xi$.
Since $\d_\xi\in E$ for $\xi<\z$ also $\d_\z\in E$.
We shall now define $\d_{\z+1}$.
We shall assume, as an induction hypothesis, that for each $\xi<\z$ the
lemma holds.
For each $\a<\d_\z$ and $k<\om$ we shall find a $\g(\a,k)$ which is
good for $(\z,\a,k)$ and we shall choose $\d_{\z+1}$ to be the least member of
$E$ greater than all these $\g(\a,k)$'s.

First we shall show that what the lemma claims holds for the case where
$\z$ is a successor or $0$.
Whenever we shall write $\z-1$ we shall assume that $\z$ is a successor.
Let $\a<\d_\z$ and $k<\om$ be given, and let
$\{\ep_n:n<\om\}=\{\ep:\ep<\d_\z\}$.
We define sequences $\enc{\a_n:n<\om}$ and $\enc{k_n:n<\om}$ so that
\item{(a)} $k_0=k$.
If $\z=0$ or $\a<\d_{\z-1}$ then $\a_0=\a$.
If $\a\ge\d_{\z-1}$ then $\a_0$ is an ordinal which is good for 
$(\z-1,\a,k)$.
In any case $\a_0<\d_\z$, \ $T_{\a_0}\su T_\a$ and
$T_{\a_0}\cap\seq k=T_\a\cap\seq k$.
\item{(b)} $\a_{n+1}<\d_\z$.
\item{(c)} $T_{\a_{n+1}}\su T_{\a_n}$.
\item{(d)} $T_{\a_{n+1}}\cap\seq{k_n}=T_{\a_n}\cap\seq{k_n}$.
\item{(e)} $T_{\a_{n+1}}$ is almost of width $f_{\ep_n}$.
\item{(f)} $k_{n+1}>k_n$ and every $\eta\in T_{\a_{n+1}}\cap\seq{k_n}$
has at least two extensions in $T_{\a_{n+1}}\cap\seq{k_{n+1}}$.

\noindent
There are indeed such sequences $\enc{\a_n:n<\om}$ and $\enc{k_n:n<\om}$.
(a) determines $k_0$ and $\a_0$; if $\a<\d_{\z-1}$ then there is an
$\a_0$ as in (a) by the induction hypothesis.
$\d_\z$ in $E$ and let us take $\a_n,\ep_n,k_n,\a_{n+1}$ for
$\a,\beta,n,\g$ in the definition of $E$, then $\d_\z\in E$ says that
there is an $\a_{n+1}$ which satisfies (b)--(e).
Since $T_{\a_{n+1}}$ is perfect there is a $k_{n+1}$ as in (f).

Let $T=\bigcap_{n\in\om}T_{\a_n}$.
By (c),(d),(f) $T$ is a perfect tree, hence it is $T_\g$ for some
$\g<\om_1$.
Since $T$, and therefore also $\gamma$, depend on $\a$ and $k$ we
denote $\g$ with $\g(\a,k)$.
As is easily seen $T_{\g(\a,k)}\su T_\a$, \
$T_{\g(\a,k)}\cap\seq k=T_\a\cap\seq k$, and for every $\ep<\d_\z$ \
$T_{\g(\a,k)}\su T_{\a_{l+1}}\su T_\a$, where $l$ is such that
$\ep=\ep_l$.
This means that (iii) of (16) holds for $\xi=\z$.
We shall have to show that (iii) holds for $\xi<\z$ and to deal
with the case where $\z$ is a limit ordinal.

If $\z$ is a limit ordinal let $\enc{\z_n:n<\om}$ be an increasing
sequence such that $\d_{\z_0}>\a$ and $\bigcup_{n<\om}\z_n=\z$.
We construct the sequences $\enc{\a_n:n<\om}$ and $\enc{k_n:n<\om}$
as in the case where $\z$ is a successor, except that (a),(b),(e) are
replaced by
\item{($\rm a'$)}$k_0=k$, \ $\a_0=\a$.
\item{($\rm b'$)}$\a_n<\d_{\z_n}$.
\item{($\rm e'$)}$\a_{n+1}$ is good for $(\z_n,\a,k)$.

\noindent
By the induction hypothesis that the lemma holds for the $\z_n$'s there
are indeed such sequences $\enc{\a_n:n<\om}$ and $\enc{k_n:n<\om}$.
Let $T=\bigcap_{n<\om}T_{\a_n}$.
As above, $T$ is a perfect tree and $T=T_{\g(\a,k)}$, \
$T_{\g(\a,k)}\su T_\a$ and $T_{\g(\a,k)}\cap\seq k=T_\a\cap\seq k$.

We shall now see that for both cases of $\z$ with which we are dealing
(iii) holds for $\xi<\z$.
If $\z$ is a successor then $\xi\le\z-1$ and since $\a_o$ is, by (a),
good for $(\z-1,\a,k)$ there is a $\beta<\d_\z$ such that
$T_{\a_0}\su T_\beta\su T_\a$ and $T_\beta$ is almost of width
$f_\ep$.
Note that if $\a<\d_{\z-1}$ then, by the induction hypothesis, we have
a $\g<\d_\z$ which is good for $(\z,\a,k)$, and if $\z=0$ then (iii)
holds vacuously, hence we may assume that $\z>0$ and 
$\a\in[\d_{\z-1},\d_\z)$.
Since $T_{\g(\a,k)}\su T_{\a_0}$ \ $\beta$ is as required by (iii).
If $\z$ is a limit ordinal then $\xi\le\z_n$ for some $n<\om$.
Since $\a_{n+1}$ is good for $\z_n$ then there is a $\beta<\d_\xi$ such
that $T_{\a_{n+1}}\su T_\beta\su T_\a$ and $T_\beta$ is almost of
width $f_\ep$.
Since $T_{\g(\a,k)}\su T_{\a_{n+1}}$ \ $\beta$ is as required by
(iii).

The only case left is that where $\zeta$ is a limit ordinal and
$\xi=\zeta$ in (iii).
Since $\a,\ep<\z$ also $\a,\ep<\z_n$ for some $n<\om$.
$\a_{n+1}$ is good for $\z_n$ hence there is a $\beta<\d_{\zeta_n}$
such that $T_{\a_{n+1}}\su T_\beta\su T_\a$ and $T_\beta$ is almost of
width $f_\ep$.
Since $T_{\g(\a,k)}\su T_{\a_{n+1}}$ and $\z_n<\z$ \ $\beta$ is as
required by (iii).
\bigbreak

\centerline{\bf Bibliography}

\noindent
[1] William Feller, An Introduction to Probability Theory and its
Applications. Wiley, New York \& London, 1950.\nl
[2] Joel Spencer, Ten Lectures on the Probabilistic Method.
CBMS-NSF Conference Series in Applied Mathematics. SIAM 1987.
\bye